\documentclass[11pt,letterpaper]{article}
\usepackage{amsmath,amssymb}
\usepackage[hidelinks]{hyperref}
\hypersetup{pdftitle={The partition sequence is not almost realizable},pdfauthor={Elias Badalov}}
\newtheorem{theorem}{Theorem}
\newtheorem{lemma}[theorem]{Lemma}

\newtheorem{proposition}[theorem]{Proposition}
\newenvironment{proof}[1][Proof]{\par\noindent\textit{#1.}\enspace}{\unskip\nobreak\hfill\(\square\)\par\medskip}
\newcommand{\F}{\mathbb F}
\title{The partition sequence is not almost realizable}
\author{Elias Badalov\\[2pt]\small Independent researcher\\[-1pt]
\small\href{mailto:ebadal.research@gmail.com}{ebadal.research@gmail.com}}
\date{November 29, 2026}
\begin{document}
\maketitle

\begin{abstract}
Geng-Rui Zhang conjectured that no positive integer multiple of the ordinary partition sequence counts the fixed points of the iterates of a single map. We prove this conjecture. For every prime \(\ell\) and every \(j\ge1\), there is an \(n\ge1\) with \(p(\ell^j n)\not\equiv p(n)\pmod\ell\). We use Euler's product and a modular-form argument modulo \(\ell\). The same obstruction applies to reciprocals of normalized integral cusp expansions and to finite products of reciprocal Euler factors.
\end{abstract}

\noindent{\small\textbf{2020 Mathematics Subject Classification.} Primary 05A17.
Secondary 11P83, 11F33, 37C25.\par}

\section{Zhang's conjecture}

A sequence \((a(n))_{n\ge1}\) is \emph{realizable} if a map \(T\) satisfies \(a(n)=\#\operatorname{Fix}(T^n)\) for every \(n\). Following Miska and Ward \cite{MiskaWard}, it is \emph{almost realizable} if a positive integer multiple is realizable. For the arithmetic of realizable sequences, see Puri and Ward \cite{PuriWard} and the survey of Byszewski, Graff, and Ward \cite{BGW}. Zhang proved that several familiar counting sequences, including the Catalan and Motzkin numbers, are not almost realizable and conjectured the same for the ordinary partition sequence \cite[Conjecture 23]{Zhang}.

\begin{theorem}[Zhang's conjecture]\label{thm:partitions}
The ordinary partition sequence \((p(n))_{n\ge1}\) is not almost realizable.
\end{theorem}

Euler's product gives an identity for the partition series modulo each prime. We use it to turn a fixed-point congruence into a statement about the reduction of a cusp form. The resulting monomial is ruled out by a consequence of B\"ocherer and Kikuta's work on modular forms modulo primes \cite{BK}, which we also prove here using Hecke operators. Sturm's bound gives an explicit upper bound for the first failed comparison \cite{Sturm}.

Ono \cite{OnoCongruences} found partition congruences on arithmetic progressions for every prime \(\ell\ge5\). Ahlgren and Boylan \cite{AB} showed that the only congruences of the form \(p(\ell n+a)\equiv0\pmod\ell\) with \(\ell\ge5\) prime and \(0\le a<\ell\) are Ramanujan's congruences modulo \(5\), \(7\), and \(11\). The fixed-point question compares values at \(n\) and \(\ell^j n\). We show that this comparison fails somewhere modulo \(\ell\) for every prime \(\ell\), including \(2\) and \(3\), and every \(j\ge1\).

\begin{lemma}[Dold congruences \cite{Dold,BGW}]\label{lem:fixed-points}
If \((a(n))_{n\ge1}\) is realizable, then for every prime \(\ell\) and \(j,n\ge1\),
\[
a(\ell^j n)\equiv a(\ell^{j-1}n)\pmod{\ell^j}.
\]
In particular, \(a(\ell^j n)\equiv a(n)\pmod\ell\).
\end{lemma}

\begin{proof}
The difference \(a(\ell^j n)-a(\ell^{j-1}n)\) counts the \(T\)-cycles whose lengths divide \(\ell^j n\) but do not divide \(\ell^{j-1}n\), with each cycle counted by its length. Every such length is divisible by \(\ell^j\). Reducing modulo \(\ell\) and applying the congruence successively gives the last assertion.
\end{proof}

\section{Proof for ordinary partitions}\label{sec:main-proof}

Write \(A(q)=\sum_{n\ge0}a(n)q^n\) for the reciprocal of a series \(E(q)\in1+q\mathbb Z[[q]]\). Let \(U_M\) take the coefficient of \(q^{Mn}\) to that of \(q^n\). Fix a prime \(\ell\) and set \(M=\ell^j\), with \(j\ge1\). Over \(\F_\ell[[q]]\), Frobenius gives
\[
A(q)=A(q^M)E(q)^{M-1}.
\]
Extracting exponents divisible by \(M\) and multiplying by \(E\) yields
\begin{equation}\label{eq:extract}
U_M(E^{M-1})-1=(U_M A-A)E.
\end{equation}
In particular \(U_M A=A\) is equivalent to \(U_M(E^{M-1})=1\) modulo \(\ell\). Multiplication by \(E\) also preserves the first nonzero index of either difference.

The following lemma is the case of B\"ocherer and Kikuta's degree-one result \cite[Proposition 4.1]{BK} that we need.

\begin{lemma}\label{lem:no-monomial}
Let \(\ell\) be prime and let \(L,w\ge1\) be integers. If \(H\in S_w(\Gamma_1(L))\) has integral Fourier coefficients at infinity, then \(H\) cannot reduce coefficientwise modulo \(\ell\) to \(c q^b\) for integers \(b\ge1\) and \(c\not\equiv0\pmod\ell\).
\end{lemma}

\begin{proof}[Hecke-operator proof]
Suppose \(H\equiv c q^b\pmod\ell\). By Dirichlet's theorem, there are infinitely many primes \(r>b\) with \(r\equiv1\pmod{L\ell}\). The diamond operator in the formula for \(T_r\) is the identity, so at \(q^m\) the two terms have indices \(rm\) and \(m/r\) in the expansion of \(H\), with the latter term taken as zero when \(m/r\) is not integral. Since \(r>b\), the first term vanishes modulo \(\ell\) at every index. Thus \(T_r H\equiv c r^{w-1}q^{rb}\pmod\ell\), and these monomials are nonzero.

The group of integral expansions in \(S_w(\Gamma_1(L))\) has finite rank. Indeed, finite dimensionality makes truncation at a sufficiently long initial block injective, so the group embeds in a finite-rank free abelian group. Its coefficientwise reduction modulo \(\ell\) is therefore finite dimensional over \(\F_\ell\) and cannot contain infinitely many linearly independent monomials.
\end{proof}

\begin{lemma}\label{lem:sieve}
If \(G(\tau)=\sum_{n\ge0}b_nq^n\) is a cusp form of weight \(w\) on \(\Gamma_1(K)\), the series obtained by keeping \(n\equiv b\pmod D\) is a cusp form of weight \(w\) on \(\Gamma_1(KD^2)\).
\end{lemma}

\begin{proof}
For \(x=t/D\), let \(R_x=\left(\begin{smallmatrix}1&x\\0&1\end{smallmatrix}\right)\).
For \(\gamma=\left(\begin{smallmatrix}a&v\\c&d\end{smallmatrix}\right)\in\Gamma_1(KD^2)\),
\[
R_x\gamma R_x^{-1}
=\begin{pmatrix}a+cx&v+(d-a)x-cx^2\\c&d-cx\end{pmatrix}
\in\Gamma_1(K).
\]
Thus \(G(\tau+x)\) is modular on the claimed group. Rational translation preserves cusp vanishing. To see this, move the source and target rational cusps to infinity. In those coordinates the change is upper triangular with positive ratio of diagonal entries, so it rescales the imaginary part by a positive constant and preserves exponential decay. Averaging the translates \(G(\tau+t/D)\) against \(e^{-2\pi i bt/D}\) over \(0\le t<D\) retains precisely the desired coefficients. The average is in characteristic zero, so integral input coefficients remain integral at infinity, including when \(\ell\mid D\).
\end{proof}

\begin{proposition}[Congruence test]\label{prop:obstruction}
Let \(E\in1+q\mathbb Z[[q]]\) and \(E^{-1}=\sum_{n\ge0}a(n)q^n\). Fix a prime \(\ell\) and \(M=\ell^j\) with \(j\ge1\). Suppose that a cusp form \(G\) of positive integral weight on \(\Gamma_1(K)\), with integral Fourier coefficients at infinity, satisfies \(G\equiv q^bE(q^d)^{M-1}\pmod\ell\) for integers \(b,d\ge1\). Then \(a(Mn)\not\equiv a(n)\pmod\ell\) for some \(n\ge1\).
\end{proposition}

\begin{proof}
Put \(D=dM\). The \emph{extracted form} \(H\) is obtained from \(G\) by retaining the terms with exponents congruent to \(b\bmod D\). By Lemma~\ref{lem:sieve}, \(H\) is cuspidal on \(\Gamma_1(KD^2)\), and its reduction is
\begin{equation}\label{eq:extracted}
H\equiv q^b\bigl(U_M(E^{M-1})\bigr)(q^D)\pmod\ell.
\end{equation}
If \(U_M A=A\) modulo \(\ell\), then \eqref{eq:extract} gives \(H\equiv q^b\pmod\ell\), contradicting Lemma~\ref{lem:no-monomial}. The constant coefficients of \(U_M A\) and \(A\) agree, so the failed congruence occurs at some \(n\ge1\).
\end{proof}

\begin{proof}[Proof of Theorem~\ref{thm:partitions}]
Euler's product gives \(E(q)=\prod_{m\ge1}(1-q^m)\) for the partition series \(A(q)=\sum_{n\ge0}p(n)q^n\). Let \(u=1\) for odd \(\ell\) and \(u=2\) for \(\ell=2\). The form
\[
\eta(24\tau)^{u(M-1)}=q^{u(M-1)}E(q^{24})^{u(M-1)}
\]
is a cusp form of positive integral weight on \(\Gamma_1(576)\). Indeed, \(u(M-1)\) is even and positive. The Newman--Gordon--Hughes eta-quotient criterion applies because \(24\) and \(576/24\) are divisible by \(24\), and Ligozat's formula gives positive orders at all cusps \cite[Theorems 1.64--1.65]{OnoBook}. See also Savitt \cite[Theorem 1.1 and Remark 1.2]{Savitt} for the criterion on \(\Gamma_1\). Its reduction is \(q^{u(M-1)}E(q^{24u})^{M-1}\), using Frobenius when \(\ell=2\). Proposition~\ref{prop:obstruction} gives \(p(\ell^j n)\not\equiv p(n)\pmod\ell\) for some \(n\ge1\).

If \(C p(n)\) were realizable for a positive integer \(C\), take a prime \(\ell\nmid C\) and use this conclusion with \(j=1\). The same congruence fails for \(C p(n)\), contradicting Lemma~\ref{lem:fixed-points}.
\end{proof}

\section{Extensions to cusp expansions and Euler products}

The same test gives a congruence obstruction for other counting sequences. In the next two results the obstruction is that for every prime \(\ell\) and integer \(j\ge1\), some \(n\ge1\) satisfies
\begin{equation}\label{eq:goal}
a(\ell^j n)\not\equiv a(n)\pmod\ell.
\end{equation}

\begin{theorem}\label{thm:general}
Let \(E(q)\in1+q\mathbb Z[[q]]\) and write \(E(q)^{-1}=\sum_{n\ge0}a(n)q^n\). Suppose that for some \(h>0\) the series \(q^h E(q)\), with \(q=e^{2\pi i\tau}\), is the expansion at infinity of a holomorphic cusp form of positive integral weight on a principal congruence subgroup \(\Gamma(N)\) for some \(N\ge1\). Then \eqref{eq:goal} holds.
\end{theorem}

A cusp form on any congruence subgroup can be restricted to a principal congruence subgroup. In particular, a form on \(\Gamma_0(N_0)\) with a Dirichlet character \(\chi\bmod N_0\) meets the group hypothesis upon restriction to \(\Gamma(N_0)\), where \(d\equiv1\pmod{N_0}\) and hence \(\chi(d)=1\).

\begin{proof}
Put \(F(\tau)=q^h E(q)\). Translation by \(N\) gives \(Nh\in\mathbb Z\). Conjugation by \(\operatorname{diag}(N,1)\) sends \(\Gamma_1(N^2)\) into \(\Gamma(N)\), so \(F(N\tau)=q^{Nh}E(q^N)\) is modular on \(\Gamma_1(N^2)\). The cusp-coordinate argument in Lemma~\ref{lem:sieve} also applies to this dilation. For \(M=\ell^j\), apply Proposition~\ref{prop:obstruction} to \(F(N\tau)^{M-1}\) with \(b=Nh(M-1)\), \(d=N\), and \(K=N^2\).
\end{proof}

\begin{theorem}\label{thm:euler}
The conclusion \eqref{eq:goal} also holds when
\[
\sum_{n\ge0}a(n)q^n=\prod_{a,m\ge1}(1-q^{am})^{-c_a},
\]
where \((c_a)_{a\ge1}\) is a finitely supported family of nonnegative integers with at least one positive entry. The sequence \((a(n))_{n\ge1}\) is therefore not almost realizable.
\end{theorem}

\begin{proof}
Write \(E(q)=\prod_{a,m\ge1}(1-q^{am})^{c_a}\). Put \(J=\sum_a a c_a>0\), let \(K_0\) be the least common multiple of the \(a\) with \(c_a>0\), and set \(M=\ell^j\). Let \(u=1\) for odd \(\ell\) and \(u=2\) for \(\ell=2\). The eta product
\[
\prod_{a\ge1}\eta(24a\tau)^{u(M-1)c_a}
=q^{uJ(M-1)}E(q^{24})^{u(M-1)}
\]
has positive integral weight and is cuspidal on \(\Gamma_1(576K_0)\). For integrality of the weight, \(M-1\) is even when \(\ell\) is odd and \(u=2\) otherwise. The two eta congruence sums vanish modulo \(24\), since \(24a\) and \(576K_0/(24a)\) are multiples of \(24\). All exponents are nonnegative and some are positive, so Ligozat's formula gives strictly positive cusp orders \cite[Theorem 1.65]{OnoBook}. The Dirichlet character supplied by the eta-quotient criterion is trivial on \(\Gamma_1(576K_0)\) \cite[Theorem 1.64]{OnoBook}\cite[Theorem 1.1 and Remark 1.2]{Savitt}.

The displayed product has reduction \(q^b E(q^{24u})^{M-1}\), where \(b=uJ(M-1)>0\). This is immediate for odd \(\ell\) and follows by squaring in characteristic two otherwise. Apply Proposition~\ref{prop:obstruction} with \(d=24u\) and \(K=576K_0\). For any positive scaling \(C\), use a prime \(\ell\nmid C\) and take \(j=1\) to contradict Lemma~\ref{lem:fixed-points} for \(C a(n)\).
\end{proof}

Theorem~\ref{thm:partitions} is the special case of Theorem~\ref{thm:euler} with \(c_1=1\) and \(c_a=0\) for \(a>1\).

\section{An effective witness bound}

\begin{proposition}
Under the hypotheses of Proposition~\ref{prop:obstruction}, let \(D=dM\) and \(L=KD^2\). Let \(H\) be the weight-\(w\) cusp form on \(\Gamma_1(L)\) obtained by retaining the terms of \(G\) with exponents congruent to \(b\bmod D\). Put
\[
B=\left\lfloor\frac{w}{12}[\operatorname{SL}_2(\mathbb Z):\Gamma_1(L)]\right\rfloor.
\]
For any prime \(r>B+b\) with \(r\equiv1\pmod{L\ell}\), there is a witness in \eqref{eq:goal} with \(n\le\lfloor(r^2-1)b/D\rfloor\).
\end{proposition}

\begin{proof}
Write \(J=\lfloor(r^2-1)b/D\rfloor\) and suppose \(a(Mn)\equiv a(n)\pmod\ell\) for every \(1\le n\le J\). The constant coefficients already agree. By \eqref{eq:extract}, \(U_M(E^{M-1})\equiv1\) through degree \(J\). Equation~\eqref{eq:extracted} then gives \(H\equiv q^b\) through exponent \(r^2 b\), since \(b+D(J+1)>r^2b\).

The coefficient at \(q^m\) of the reduction of \(T_r H\) is the sum of the original coefficient at \(q^{rm}\) and \(r^{w-1}\) times that at \(q^{m/r}\), with the latter zero for nonintegral indices. For \(m\le B\) both terms vanish because \(r>B+b\). At \(m=rb\) the coefficient is \(r^{w-1}\not\equiv0\pmod\ell\). This contradicts Sturm's bound at level \(\Gamma_1(L)\) \cite{Sturm}.
\end{proof}

\paragraph{Lean verification.}
The arXiv ancillary archive \texttt{zhang-conjecture-23-lean.zip} contains a Lean~4 formalization of Theorem~\ref{thm:partitions} using prime-square dilation. It includes the source code, pinned dependencies, build instructions, and a record of the kernel's axiom check. The formalization covers Zhang's conjecture. The assertion for every dilation exponent, the explicit witness bound, and the generalizations in this paper are proved in the text.

\paragraph{AI assistance.} The proofs in this paper were developed and checked with assistance from AI tools provided by OpenAI, Anthropic, and Google.

\end{document}